# PARAMETRIC VIBRATIONS OF A DAMAGED ORTHOTROPIC GEOMETRICAL SHELL STIFFENED WITH AN INHOMOGENEOUS ROD AND RINGS ON VISCOELASTIC MEDIUM

I.G. Aliyev [1]   F.S. Latifov [1]   A.M. Guliyeva [1]   N.H. Seyidov [1]   A.I. Hasanov [2]

*1. Azerbaijan Architectural and Construction University, Baku, Azerbaijan*
*i_q_aliyev@mail.ru, flatifov@mail.ru, azizaquliyeva1960@gmail.com, kristal_namiq@mail.ru*
*2. Nakhichevan State University, Nakhichevan, Azerbaijan, hesen.1955@mail.ru*

**Abstract.** The paper considers a problem of parametric vibrations of a viscoelastic medium-contacting damaged cylindrical shell under the action of external force $p = p_0 + p_1 \sin \omega_1 t$ ( $p_0$ is a medium or main force, $p_1$ is the amplitude of force change, $\omega_1$ is the frequency of change of the variable part of the force) and stiffened along its generatrix with a rod whose material is inhomogeneous and in circular direction with inhomogeous rings. The structural element considered in the presented article consists, according to the geometric structure, of a coating and reinforcement elements, according to the mechanical characteristics of heterogeneous coatings along the length, having damage inside due to their physical structure and, finally, a system in contact with a viscoelastic medium. Taking into account one of the existing models (Winkler or Pasternak), the contact conditions between the coating and the reinforcement elements and the influence of the medium on the coating, the frequency equation of oscillation was solved, the results were analyzed. Theory of hereditary type damage under the action of external force is often used for taking into account the damages in the structure of a cylindrical shell forced to vibration. The Hamilton-Ostrogradsky variational principle is used for solving the problem. The results obtained can be used in the foundations of bridges built across mountain rivers. Such support is cost-effective. It should be noted that reinforcement and concrete mortar are also used in the internal parts of the cylindrical supports used. The article proposes to fill the inside of the cylindrical lid under study with clay.

## 1. INTRODUCTION

Analysis of medium contacting constructions taking into account external factors in working condition is of great importance in practice. Ensuring the serviceability of shell construction is an urgent issue. Long-time exposure of the construction to vibration causes damage and accumulation of damages in the internal area. Therefore, the study of parametric vibrations in viscoelastic medium taking into amount these factors and the inhomogeneity of the material of the structure is one of the urgent problems.

In [1-3], a problem on parametric vibrations of a viscoelastic medium contacting cylindrical shell subjected to the action of external forces and stiffened with homogenous ribs is considered. Three cases of location of the rods on the surface of the cylindrical shell were considered 1) stiffening with rods 2) stiffening by means of rings 3) In this case, the anchoring elements form an orthogonal grid. The model was modeled in the viscoelastic form and was studied by means of the Lame equation in elastic effect displacements.

A frequency equation was constructed to find the frequency of parametric oscillations of the system using the contact condition, and an asymptotic approach. Optimization parameter was introduced and optimal variant of the number of ribs was found [4]. Studies parametric vibrations of damaged structural elements made of isotropic plane material and modeled as a rod or shell. [5] the article studies parametric oscillations of a reinforced anisotropic coating in a geometric nonlinear formulation, which were solved by the energy method.

The dependence of the ratio of nonlinear frequency to linear frequency on the curvature of the shell was determined for different number of rods. The paper [6] considers vibrations of a homogeneous, damaged, stiffened, viscoelastic, medium, filled cylindrical shell. The paper [16] considers vibrations of an amplified, damaged anisotropic cylindrical shell under axial compression by a flowing liquid. The carried-out analysis shows that this problem has not found its solution for a damaged cylindrical shell with orthotropic features, stiffened with inhomogeneous rods and rings and dynamically contacting with viscoelastic medium. Calculations show that it is important to take into account the adequate properties of the elements that make up the structure, including the heterogeneity of the shafts and rings used in the reinforcement, the orthotropy of the material of the cylindrical shell, and viscoelasticity.

The inner area is filled with clay, and the cylindrical coating is continuously exposed to external forces, leading to potential damage. The found crisis force makes it possible to control the performance of the structure.

## 2. FORMULATION OF THE PROBLEM

For studying the is a problem in which is considered in this article under the action of a variable external force and stiffened with inhomogeneous rings and rods, we will use the Hamilton, Ostrogradskiy variational principle (Figure 1). Hereditary type damage theory was used for taking into account damages. According to this theory, strain components in homogeneous body are determined as follows [7]:

$$\varepsilon_{ij} = \overline{\varepsilon}_{ij} + M^* \times \sigma_{ij} \quad (1)$$

where, $M^*$ is a hereditary type integral operator describing the damage process and is as follows:

$$M^* \times \sigma_{ij} = \sum_{k=1}^{n} f(t_k^+) \int_{t_k^-}^{t_k^+} M(\overline{x}, t_k^+ - \tau)\sigma_{ij}(\tau)d\tau + \int_{t_{n+1}^-}^{t} M(\overline{x}, t_k^+ - \tau)\sigma_{ij}(\tau)d\tau \quad (2)$$

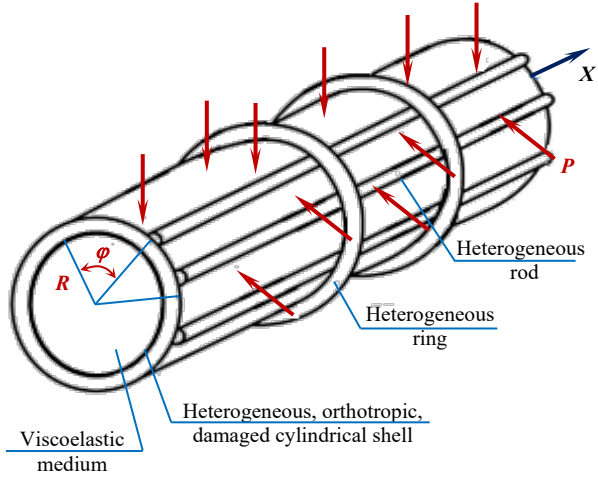

Figure 1. Physical model of a damaged orthotropic cylindrical shell

In the expression (2) $M(\overline{x}, t_k^+ - \tau)$ is a damage kernel, $(t_k^-, t_k^+)$ is a time interval during which the active stress insures an increase in damage, $f(t_k^+)$ is a recovery function and depends on the volume of damages accumulated in a cycle. The value $f(t_k^+) = 0$ of this function corresponds to complete recovery of damages, the value $f(t_k^+) = 1$ to the absence of recovery process. The values between zero and unit express partial recovery of damages. To identify the interval $(t_k^-, t_k^+)$ we need to give special condition. This condition depends on specific features of the structure, its working conditions and the type of loading. Taking into account the Equation (2), for the total energy of a cylindrical shell with damages taken into account we can write:

$$J = \frac{1}{2}R^2 \int_{x_1}^{x_2}\int_{y_1}^{y_2} \{N_{11}\varepsilon_{11} + N_{22}\varepsilon_{22} + N_{12}\varepsilon_{12} - M_{11}\chi_{11} - M_{22}\chi_{22} - M_{12}\varepsilon_{12} +$$

$$+ N_{11}(\sum_{k=0}^{n} f(t_k^+) \times \int_{t_k^-}^{t_k^+} M(\overline{x}, t_k^+ - \tau)N_{11}d\tau + \int_{t_{n+1}^-}^{t} M(\overline{x}, t - \tau)N_{11}d\tau) +$$

$$+ N_{22}\left(\sum_{k=0}^{n} f(t_k^+) \int_{t_k^-}^{t_k^+} M(\overline{x}, t_k^+ - \tau)N_{22}d\tau + \int_{t_{n+1}^-}^{t} M(\overline{x}, t - \tau)N_{22}d\tau\right) +$$

$$+ N_{12}\left(\sum_{k=0}^{n} f(t_k^+) \times \int_{t_k^-}^{t_k^+} M(\overline{x}, t_k^+ - \tau)N_{12}d\tau + \int_{t_{n+1}^-}^{t} M(\overline{x}, t - \tau)N_{12}d\tau\right) -$$

$$- M_{11}\left(\sum_{k=0}^{n} f(t_k^+) \int_{t_k^-}^{t_k^+} M(\overline{x}, t_k^+ - \tau)M_{11}d\tau + \int_{t_{n+1}^-}^{t} M(\overline{x}, t - \tau)M_{11}d\tau\right) -$$

$$- M_{22}\left(\sum_{k=0}^{n} f(t_k^+) \times \int_{t_k^-}^{t_k^+} M(\overline{x}, t_k^+ - s\tau)M_{12}d\tau + \int_{t_{n+1}^-}^{t} M(\overline{x}, t - \tau)M_{12}d\tau\right) -$$

$$- M_{12}\left(\sum_{k=0}^{n} f(t_k^+) \int_{t_k^-}^{t_k^+} M(\overline{x}, t_k^+ - \tau)M_{12}d\tau + \int_{t_{n+1}^-}^{t} M(\overline{x}, t - \tau) \times\right.$$

$$\left. \times M_{12}d\tau\right) + \rho_0 h \int_{x_1}^{x_2}\int_{y_1}^{y_2}\left[\left(\frac{\partial u}{\partial t}\right)^2 + \left(\frac{\partial \vartheta}{\partial t}\right)^2 + \left(\frac{\partial w}{\partial t}\right)^2\right]dxdy \quad (3)$$

The internal forces and moments included in Equation (3) will be taken in the following form [13]:

$$N_{ij} = \int_{-h/2}^{h/2}(\sigma_{ij} + zw_{ij})dz, \quad M_{ij} = -\int_{-\frac{h}{2}}^{\frac{h}{2}}(\sigma_{ij} + zw_{ij})zdz \quad (4)$$

where, $w_{11}, w_{22}, w_{21} = w_{12}$ and $\sigma_{11}, \sigma_{22}, \sigma_{12}, \varepsilon_{11}, \varepsilon_{22}, \varepsilon_{12}$ formulas were obtained to calculate the Equation [3]. Let's express the internal forces, stresses and displacements included in Equation (4) by deformations as follows:

$$\chi_{11} = \frac{\partial^2 w}{\partial x^2}; \quad \chi_{22} = \frac{\partial^2 w}{\partial y^2}; \quad \chi_{12} = -2\frac{\partial^2 w}{\partial x \partial y}; \quad b_{11} = \frac{E_1}{1 - v_1 v_2}$$

$$b_{66} = G; \quad b_{22} = \frac{E_2}{1 - v_1 v_2}; \quad b_{12} = \frac{v_2 E_1}{1 - v_1 v_2} = \frac{v_1 E_2}{1 - v_1 v_2} \quad (5)$$

For taking into account the inhomogeneity of the rods, we will consider trat the modules of elasticity and density of the material is a function of the coordinate $x$. In this case, for the total energy of the rods stiffened we can write:

$$J_i = \frac{1}{2}\sum_{i=1}^{k_1}\int_0^1\left[\tilde{E}_i(x)F_i\left(\frac{\partial u_i}{\partial x}\right)^2 + \tilde{E}_i(x)J_{yi}\left(\frac{\partial^2 w_i}{\partial x^2}\right)^2 +\right.$$

$$\left. + \tilde{E}_i(x)J_{zi}\left(\frac{\partial^2 \vartheta_i}{\partial x^2}\right)^2 + \tilde{G}_i(x)J_{kpi}\left(\frac{\partial \varphi_{kpi}}{\partial x}\right)^2\right]dx +$$

$$+ \sum_{i=1}^{k_1}\rho_i(x)F_i\int_0^1\left[\left(\frac{\partial u_i}{\partial t}\right)^2 + \left(\frac{\partial \vartheta_i}{\partial t}\right)^2 + \left(\frac{\partial w_i}{\partial t}\right)^2 + \frac{J_{kpi}}{F_i}\left(\frac{\partial \varphi_{kpi}}{\partial t}\right)^2\right]dx \quad (6)$$

For taking into account the inhomogeneity of the rings, we will consider the modules of elasticity and density of the material as a function of the coordinate $\varphi$. In this case for the total energy of rings stiffened in the surface perpendicular to axis of the cylindrical shell, we can write:

$$J_j = \frac{R}{2}\sum_{j=1}^{k_2}\int_0^{2\pi}\left[\bar{E}_j(\varphi)F_j\left(\frac{\partial\vartheta_j}{r\partial\varphi}-\frac{w_j}{r}\right)^2 + \bar{E}_j(\varphi)J_{zj}\times\right.$$

$$\times\left(\frac{\partial^2 w_j}{\partial x^2}+\frac{w_j}{R^2}\right)^2 + \bar{E}_j(\varphi)J_{zj}\left(\frac{\partial^2 u_j}{R^2\partial\varphi^2}-\frac{\varphi_{kpj}}{R}\right)^2 +$$

$$\left.+\bar{G}_j(\varphi)J_{kpj}\left(\frac{\partial\varphi_{kpj}}{R\partial\varphi}+\frac{1}{R}\frac{\partial u_j}{\partial\varphi}\right)^2\right]d\varphi + \quad (7)$$

$$+R\sum_{j=1}^{k_2}F_j\int_0^{2\pi}\bar{\rho}_j(\varphi)\left[\left(\frac{\partial u_j}{\partial t}\right)^2 + \left(\frac{\partial\vartheta_j}{\partial t}\right)^2 + \left(\frac{\partial w_j}{\partial t}\right)^2 +\right.$$

$$\left.+\frac{J_{kpj}}{F_j}\left(\frac{\partial\varphi_{kpj}}{\partial t}\right)^2\right]d\varphi$$

An explanation of all the notations included in Equation (3)-(7) are listed in the articles [3-6], [14], [15]. In this article, these layouts are saved as is. The work $A_0$ performed from the impact of a viscoelastic medium from the inside on the cylindrical coating, the work $A_p$ performed from the impact of an intense load from the outside, respectively, will be calculated as follows:

$$A_0 = -R\int_0^l\int_0^{2\pi}q_z w\, dx\, d\varphi \quad (8)$$

$$A_p = -4R\int_0^l\int_0^{\pi/4}pw\, dx\, d\varphi \quad (9)$$

The force contained in the Equation (9) and acting as viewed from the medium to cylindrical shell is determined:

$$q_z = k_\vartheta w - k_p\left(\frac{\partial^2 w}{\partial x^2}+\frac{\partial^2 w}{R\partial\varphi^2}\right) - \int_0^t \Gamma(t-\tau)w(\tau)d\tau$$

where, $k_\vartheta$ is a Winkler coefficient, $k_p$ is a Pasternak coefficient and are found experimentally, $t$ is time, $\Gamma(t-\tau)$ is a viscosity kernel. It is considered that rods and rings were rigidly hinged to the cylindrical shell. In this case the following contact conditions are fulfilled [3].

### 3. SOLUTION OF THE PROBLEM

For $x=0$ and $x=1$, within the boundary conditions $\vartheta = w = T_1 = M_1 = 0$ we will look for shell displacements:

$$u = u_0 \cos n\varphi \cos mx \sin\omega t$$
$$\vartheta = \vartheta_0 \sin n\varphi \sin mx \sin\omega t \quad (10)$$
$$w = w_0 \cos n\varphi \sin mx \sin\omega t$$

where, $u_0, \vartheta_0, w_0$ are unknown constants. By means of Equation (10), the active loading period included in the damage operator can be determined from the decreasing condition of these functions:

$$\left[\left(\frac{\pi}{2}+2\pi k\right)/\omega; \left(\frac{3\pi}{2}+2\pi k\right)/\omega\right]$$

The characteristically time $T$ is determined as the greatest of the times $t_n^+$. Taking into account the Equation (10) in the Equations (3), (6), (7), (8) and (9), for the Hamilton action $W = \int_0^T J dt$ we obtain $M(\bar{x}, t-\tau) = \gamma = \text{const.}$

$$W = \frac{\pi LhR^2}{4}\left\{\left[\left(k^2 b_{11}+\frac{n^2}{R^2}b_{66}\right)\left(\frac{T}{2}-\frac{\sin 2\omega T}{4\omega}\right)-\right.\right.$$

$$-\frac{\gamma h}{\omega}F(T)T_1 + \frac{\gamma h^3}{16\omega}F(T)T_1\right]u_0^2 + \left[\left(\frac{n^2}{R^2}b_{22}+k^2 b_{66}\right)+\right.$$

$$+\left(\frac{T}{2}-\frac{\sin 2\omega T}{4\omega}\right)-\frac{\gamma h}{\omega}F(T)T_2 + \frac{\gamma h^3}{16\omega}F(T)T_2\right]\vartheta_0^2 +$$

$$+\left[\left(b_{22}-\frac{hk^2 b_{12}}{4}-\frac{hn^2 b_{22}}{4R^2}-\frac{hk^2}{2}\times\left(\frac{hb_{11}k^2}{3}+\frac{hn^2 b_{12}}{3R^2}-b_{12}\right)-\right.\right.$$

$$-\frac{h}{2}\left(\frac{hn^2 k^2 b_{12}}{3R^2}+\frac{hn^4 b_{22}}{R^4}-\frac{n^2 b_{22}}{R^2}\right)-\frac{h^2 n^2 k^2 b_{66}}{6R^2}\right)\times$$

$$\times\left(\frac{T}{2}-\frac{\sin 2\omega T}{4\omega}\right)-\frac{\gamma h F(t)}{\omega}\left(T_3-T_4-\frac{h^3 n^2 k^2 b_{66}}{4R^2}\right)+$$

$$+\frac{\gamma h^3}{16\omega}F(T)\left(T_3+T_4+\frac{4h^3 n^2 k^2 b_{66}}{9R^2}\right)\right]w_0^2 +$$

$$+\left[\frac{2nk}{R}(-b_{12}-b_{66})\left(\frac{T}{2}-\frac{\sin 2\omega T}{4\omega}\right)+\frac{2nk}{R}\frac{\gamma h^2}{\omega}F(T)\times\right.$$

$$\times(b_{11}b_{12}+b_{11}b_{22})-\frac{\gamma h^4}{16\omega}F(T)\frac{2nk}{R}b_{11}b_{12}+b_{12}b_{22}+$$

$$+b_{66}^2\right]u_0\vartheta_0 + \left[\left(-2kb_{12}+\frac{hk^3 b_{11}}{4}+\frac{hn^2 kb_{12}}{4R^2}-\frac{hn^2 kb_{66}}{2R^2}-\right.\right.$$

$$-\frac{kh^2}{2\omega}\left(b_{11}k^2+\frac{n^2 b_{12}}{R^2}-\frac{2n^2 b_{66}}{R^2}\right)\right)\left(\frac{T}{2}-\frac{\sin 2\omega T}{4\omega}\right)+$$

$$+\left[\left(\frac{hk^2 nb_{66}}{2R}+\frac{hk^2 nb_{12}}{4R}+\frac{h^2 n^3}{4R^3}b_{22}-\right.\right.$$

$$-\frac{khn}{R}b_{66}+\frac{2b_{22}n}{R}\right)\left(\frac{T}{2}-\frac{\sin 2\omega T}{4\omega}\right)-$$

$$-\frac{2n\gamma h^2}{R\omega}F(T)\left(b_{12}T_3+b_{22}T_4+\frac{hk^2 b_{66}^2}{2R}\right)+$$

$$+\frac{n\gamma h^4}{8R\omega}F(T)(b_{12}T_5+b_{22}T_6)\right]u_0 w_0\right\}+$$

$$+\rho_0 h\frac{\pi\omega^2 L}{2}(u_0^2+\vartheta_0^2+w_0^2)\left(\frac{T}{2}+\frac{\sin 2\omega T}{4\omega}\right)-$$

$$-\frac{4}{nk}(\cos kL-1)\sin\frac{n\pi}{4}\left[-\frac{p_0}{\omega}(\cos\omega T-1)+\right.$$

$$+\frac{1}{2}p_1\left(\frac{2}{\omega-\omega_1}\sin(\omega-\omega_1)T-\frac{2}{\omega+\omega_1}\times\right.$$

$$\times\sin(\omega+\omega_1)T\right]w_0 + \left\{\frac{m^2\pi^2}{2l^2}\sum_{i=1}^{k_1}[F_i I_1 \sin^2 n\varphi_i u_0^2 +\right.$$

$$+(J_{xi}I_2+J_{kpi}I_3)\cos^2 n\varphi_i \vartheta_0^2 + (J_{zi}I_2+J_{kpi}I_3)\times$$

$$\times\sin^2 n\varphi_i w_0^2 + kJ_{kpi}I_3 \sin 2n\varphi_i \vartheta_0 w_0]+$$

$$+\omega^2 \sum_{i=1}^{k_1} F_i \left[ I_4 \sin^2 kn u_0^2 + I_5 \left(1 + \frac{J_{kpi}}{F_i R^2}\right) \cos^2 n\varphi_i \vartheta_0^2 + \right.$$

$$+ I_5 \left(1 + \frac{J_{kpi} k^2}{F_i R^2}\right) \sin^2 n\varphi_i w_0^2 +$$

$$\left. + I_5 \frac{J_{kpi}}{F_i R^2} \sin 2n\varphi_i \vartheta_0 w_0 \right] \left(\frac{T}{2} - \frac{\sin 2T}{4\omega}\right) -$$

$$- 2\pi R \left\{ \frac{\omega l}{4(\omega^2 + \psi^2)} \left[ \frac{1}{\omega}(1 - \cos\omega T) + \frac{1}{\psi}(e^{\psi T}\sin^2 \omega T - \right.\right.$$

$$\left.\left. - \frac{\omega}{1+\psi^2}(\psi - \frac{1}{2\omega}(e^{\psi T}\cos 2\omega T - 1))\right] w_0^2 + \right.$$

$$+ \sum_{j=1}^{k_2} \left\{ \frac{R \cos^2 mx_j}{2}(J_{zj} I_6 + \frac{n^2}{R^2} J_{kp.j} I_9 + \right.$$

$$+ R\omega^2 F_j I_7) u_0^2 + \left(\frac{n^2}{2R} I_6 + \frac{R^2}{2}\omega^2 F_j I_8\right) \sin^2 mx_j \vartheta_0^2 +$$

$$+ \left[\left(\frac{1}{2R} + \frac{R}{2} J_{zj}\left(\frac{1}{R^2} - m^2\right)\right) I_6 \sin^2 mx_j + \right.$$

$$+ \frac{1}{2R}\left(J_{zj} m^2 I_6 + n^2 m^2 J_{kp.j} I_9\right)\cos^2 mx_j +$$

$$+ R\omega^2 F_j \left(\sin^2 mx_j + \frac{J_{kp.j}}{F_j}\cos^2 mx_j\right) I_7 \right] w_0^2 - \frac{1}{R} \times$$

$$+ R\omega^2 F_j \left(\sin^2 mx_j + \frac{J_{kp.j}}{F_j}\cos^2 mx_j\right) I_7 \right] w_0^2 - \frac{1}{R} \times$$

$$\times mn^2 \left(\frac{1}{R} J_{zj} I_6 - J_{kp.j} I_9\right)\cos^2 mx_j u_0 w_0 -$$

$$- \frac{n}{R} I_6 \sin^2 mx_j \vartheta_0 w_0 \right\} \left(\frac{T}{2} - \frac{\sin 2T}{4\omega}\right)$$

where,

$$I_1 = \int_o^l \widetilde{E}_i(x)\cos^2 mx\, dx; \quad I_2 = \int_o^l \widetilde{E}_i(x)\sin^2 mx\, dx;$$

$$I_3 = \int_o^l \widetilde{G}_i(x)\sin^2 mx\, dx; \quad I_4 = \int_0^l \widetilde{\rho}_i(x)\cos^2 mx\, dx;$$

$$I_5 = \int_0^l \widetilde{\rho}_i(x)\sin^2 mx\, dx; \quad I_6 = \int_o^{2\pi} \overline{E}_j(\varphi)\cos^2 n\varphi\, d\varphi;$$

$$I_7 = \int_o^{2\pi} \rho_j(\varphi)\cos^2 n\varphi\, d\varphi; \quad I_8 = \int_o^{2\pi} \rho_j(\varphi)\sin^2 n\varphi\, d\varphi;$$

$$I_9 = \int_o^{2\pi} \overline{G}_j(\varphi)\sin^2 n\varphi\, d\varphi$$

$$F(T) = \frac{1}{2\omega}\left(\sin^2 \omega T + 4R_t \sin^2 \frac{\omega T}{2}\right)$$

It should be noted that in Equation (11), taking the terms with indices $i$ equal to zero, we obtain a problem of parametric vibration of a damaged, viscoelastic medium-contacting cylindrical shell stiffened with inhomogeneous rings only in the circular direction, taking the terms with indices $j$ equal to zero we obtain a problem of parametric vibrations of a viscoelastic, medium-contacting, damaged cylindrical shell stiffened along the generatrix with rods whose material is inhomogeneous [12]. Using the stationing conditions of the Ostrogrdasky, Hamilton action $\delta W = 0$ with respect to the constants $u_0$, $\vartheta_0$, $w_0$ we obtain the following system of inhomogeneous equation:

$$\begin{cases} l_{11}u_0 + l_{12}\vartheta_0 + l_{13}w_0 = 0 \\ l_{21}u_0 + l_{22}\vartheta_0 + l_{23}w_0 = 0 \\ l_{31}u_0 + l_{32}\vartheta_0 + l_{33}w_0 = \breve{\varphi}_* \end{cases} \quad (12)$$

where,

$$\breve{\varphi}_* = -\frac{4}{nk}(\cos kL - 1)\sin\frac{n\pi}{4}\left[-\frac{p_0}{\omega}(\cos\omega T - 1) + \right.$$

$$\left. + \frac{1}{2}p_1\left(\frac{2}{\omega-\omega_1}\sin(\omega-\omega_1)T - \frac{2}{\omega+\omega_1}\sin(\omega+\omega_1)T\right)\right]$$

Having solved the system (13), we can determine the constants $u_0$, $\vartheta_0$, $w_0$.

$$u_0 = \frac{\Delta_1}{\Delta}\ ; \quad \vartheta_0 = \frac{\Delta_2}{\Delta}\ ; \quad w_0 = \frac{\Delta_3}{\Delta} \quad (13)$$

In the expressions (13), $\Delta$ is the principal determinant of the system (12), $\Delta_i (i = 1,2,3,)$ are auxiliary determinants. Having substitute the expressions (13) in (10), for the displacement of the points of the shell we obtain [13]:

$$u = \frac{\Delta_1}{\Delta}\cos n\varphi \cos mx \sin \omega t$$

$$\vartheta = \frac{\Delta_2}{\Delta}\sin n\varphi \sin mx \sin \omega t \quad (14)$$

$$w = \frac{\Delta_3}{\Delta}\cos n\varphi \sin mx \sin \omega t$$

To determine the vertical force, we use the equality $\frac{\Delta_3}{\Delta} = w_0$. From this equality we obtain [14]:

$$\breve{\varphi}_* = \frac{\Delta w_0}{l_{11}l_{22} - l_{21}l_{12}} \quad (15)$$

Taking into account the expression (12) of $\breve{\varphi}_*$ in (15), we can find the force $p_1$:

$$p_1 = \frac{\Delta w_0}{\alpha_{22}(l_{11}l_{22} - l_{21}l_{12})} - \frac{\alpha_{11}p_0}{\alpha_{22}} \quad (16)$$

where,

$$\alpha_{11} = \frac{4}{nk\omega}(\cos kL - 1)\sin\frac{n\pi}{4}(\cos\omega T - 1)$$

$$\alpha_{22} = -\frac{2}{nk}(\frac{2}{\omega-\omega_1}\sin(\omega-\omega_1)T -$$

$$- \frac{2}{\omega+\omega_1}\sin(\omega+\omega_1)T)(\cos kL - 1)\sin\frac{n\pi}{4}$$

## 4. NUMERICAL RESULTS

Changing the wave numbers $n$ and $m$, we calculate $p_1$ and choosing $(p_1)_{\min}$ from them we find the critical force $p_{1b}$. When there are no damages, in the expression $F(T) = \frac{1}{2\omega}\left(\sin^2 \omega T + 4R_t \sin^2 \frac{\omega T}{2}\right)$ taking $l_t, R_t = 0$, we can calculate the critical value of the force $p_{1b}$. The expression (16) was numerically calculated for the following values of variables:

$\tilde{E}_{j0} = 6.67 \times 10^9 \text{ N/m}^2$ , $h = 0.45 \text{ mm}$ , $v_1 = 0.11$

$v_2 = 0.19$ , $l = 800 \text{ mm}$ , $R = 160 \text{ mm}$ , $\beta = 0.05$

$A = 0.1615$ , $\rho_0 = \rho_j = 7.8 \text{ g/cm}^3$ , $\omega_1 = 2\omega$

$h_i = 1.39 \text{ mm}$ , $w_0 = 0.1 \text{ mm}$ , $\omega = 100 \text{ hs}$

$I_{kp.i} = 0.23 \text{ mm}^4$ , $I_{zi} = 1.3 \text{ mm}^4$ , $F_i = 5.2 \text{ mm}^2$

$k_\vartheta = \dfrac{10^6 \text{ N}}{\text{m}^3}$ , $k_p = \dfrac{10^4 \text{ N}}{\text{m}}$ , $h_j = 1.39 \text{ mm}$

$I_{kp.j} = 0.48 \text{ mm}^4$ , $I_{xj} = 19.9 \text{ mm}^4$ , $G = 3.5 \dfrac{10 \times N}{\text{m}^2}$

It was accepted that

$\tilde{E}_j(\varphi) = \tilde{E}_{j0}\left(1 + \sigma \dfrac{\varphi}{2\pi}\right)$ , $\tilde{G}_j(\varphi) = \tilde{G}_{j0}(1 + \mu \dfrac{\varphi}{2\pi})$

$\tilde{\rho}_j(\varphi) = \tilde{\rho}_{j0}\left(1 + \tau \dfrac{\varphi}{2\pi}\right)$ , $\sigma = \mu = \tau = 0.4$

where, $\sigma, \mu, \tau$ are inhomogeneity parameters.

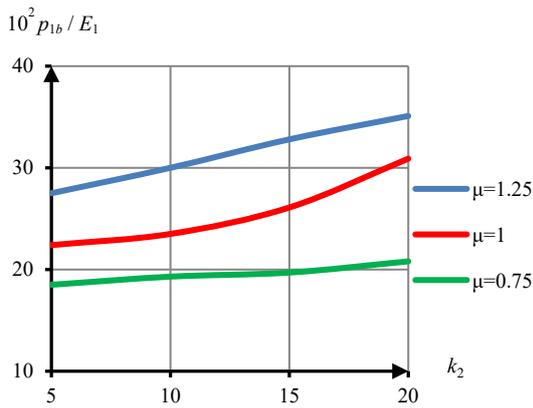

Figure 2. Dependence $P_{1b}/E_1$ on the number of rings $k_2$

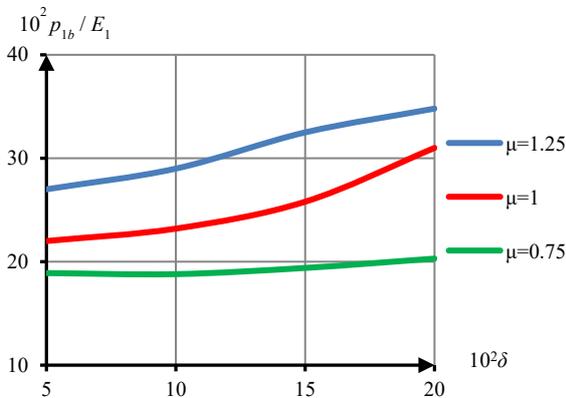

Figure 3. $10^2 p_{1b}/E_1$ and the relationship between

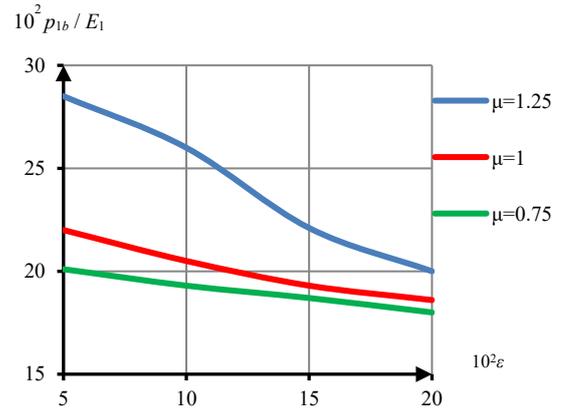

Figure 4. The relationship between $p_{1b}/E_1$ and $\varepsilon$

The result of calculations were given in Figures 2-4 in the form of dependence of the ratio $p_{1b}/E_1$ on the number of the rings $k_2$ of the inhomogeneity parameters $\sigma$, $\tau$ and $\mu = E_1/E_2$. As can be seen from Figure 2, as the number of rings increases, the value of the critical force increases.

The article presents a parallel comparison of the critical force found for the stability problem with the frequency parameter, which is the main indicator in oscillation problems. That is, since an increase in the number of reinforcement elements in the issue of stability increases rigidity, the initial critical force increases accordingly. In some cases, this approach is accompanied by contraindications. In oscillatory problems, an increase in heterogeneity due to an increase in the number of reinforcement elements creates an inertia effect, as a result of which the oscillation frequency parameter decreases. Note that in a special case, when there no rings in stiffening, the results coincide with result of the paper [11].

Taking into account the heterogeneity of the reinforced elements allows you to save a sufficient amount of material while optimizing the geometric dimensions of the structure under study, as well as its weight. Such a comparison was carried out and analyzed with the works [6, 15], where due to the heterogeneity of the reinforcement elements, material consumption decreased by 6-9%, and the structure became lighter while maintaining strength.

## NOMENCLATURES

**Symbols / Parameters**

$J$ : The energy of the body of the structure;

$J_i$ , $J_j$ : Accordingly, energy of reinforcement elements;

$E_i(x)$ , $F_i$ , $E_j(\varphi), F_i$ : Geometric characteristics of reinforcement elements;

$u_i, \vartheta_i, w_i$ , $u_j, \vartheta_j, w_j$ : They are the displacements of the points of the reinforcement elements;

$\rho_j(\varphi), G_i, J_{jkp}, J_{xj}$ : The density, modulus of elasticity and moments of inertia of the transverse reinforcement element during torsion;

$u, \vartheta, w$ and $R, h$ : Displacements and geometric dimensions of an arbitrary point of body of the structure.

## BIOGRAPHIES

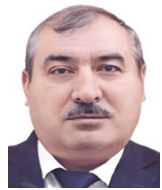

Name: **Ilgar**
Middle Name: **Giyas**
Surname: **Aliyev**
Birthday: 15.05.1966
Birthplace: Amasia, Armenia
Master: Heating Gas Supply and Ventilation, Faculty of Water Management and Engineering Communication Systems, Azerbaijan University of Architecture and Construction, Baku, Azerbaijan, 1988
Doctorate: Heating Gas Supply and Ventilation, Construction and Operation of Oil and Gas Pipelines Department, Water Management and Engineering Communication Systems Faculty, Azerbaijan Architecture and Construction University, Baku, Azerbaijan, 1994
The Last Scientific Position: Assoc. Prof., Faculty of Civil Engineering, Azerbaijan University of Architecture and Construction, Baku, Azerbaijan, Since 1997
Research Interests: Construction and Operation of Oil-Gas Pipelines, Operation and Reconstruction of Buildings and Facilities, Designing Buildings and Facilities
Scientific Publications: 46 Papers, 11 Books, 1 Patent, 12 Projects, 15 Theses

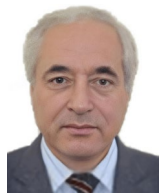

Name: **Fuad**
Middle Name: **Seyfaddin**
Surname: **Latifov**
Birthday: 01.07.1955
Birthplace: Baku, Azerbaijan
Master: Mechanics, Mechanical and Mathematical Faculty, Baku State University, Baku, Azerbaijan, 1977
Doctorate: Ph.D., Department of Mechanics, Leningrad State University, St. Petersburg, Russia, 1983
Doctorate: Dr. Sci., Mathematics, Deformed Solid Mechanics, Department of Mechanics, Institute of Mathematics and Mechanics, Azerbaijan National Academy of Sciences, Baku, Azerbaijan, 2003
The Last Scientific Position: Prof., Department of Mechanics, Faculty of Water Management and Engineering Communication Systems, Azerbaijan University of Architectue and Construction, Baku, Azerbaijan, Since 2012

Research Interests: Mechanics, Mathematics
Scientific Publications: 95 Papers, 12 Books, 30 Theses
Scientific Memberships: Member of the Turkish World Research International Academy of Sciences, Laureate of International Gold Star Medal

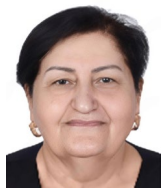
Name: **Aziza**
Middle Name: **Mursal**
Surname: **Guliyeva**
Birthday: 02.09.1960
Birthplace: Saatli, Azerbaijan
Master: Faculty of Construction, Azerbaijan University of Architecture and Construction, Baku, Azerbaijan, 1982
Doctorate: Ph.D., Building Structures, Azerbaijan University of Architecture and Construction, Baku, Azerbaijan, 1987
The Last Scientific Position: Assoc. Prof., Department of Reinforced Concrete Structures, Azerbaijan University of Architecture and Construction, Baku, Azerbaijan, Since 2000
Scientific Interests: Reinforced Concrete, Metal and Wooden Structures, Structural Mechanics Scientific, Metal Science
Scientific Publications: 45 Papers, 3 Textbooks

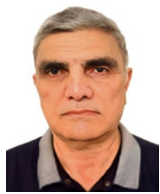
Name: **Namiq**
Middle Name: **Hasan**
Surname: **Seyidov**
Birthday: 27.01.1958
Birthplace: Nakhichevan, Azerbaijan
Master: Faculty of Construction, Azerbaijan University of Architecture and Construction, Baku, Azerbaijan, 1975
Doctorate: Ph.D., Building Structures, Moscow Engineering and Construction Institute, Moscow, Russia, 1988
The Last Scientific Position: Assoc. Prof., Department of Reinforced Concrete Structures, Azerbaijan University of Architecture and Construction, Baku, Azerbaijan, Since 1998
Research Interests: Reinforced Concrete, Metal and Wooden Structures, Structural Mechanics
Scientific Publications: 70 Papers, 1 Textbook

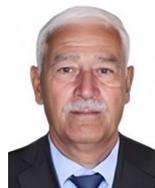
Name: **Abdulla**
Middle Name: **Ismail**
Surname: **Hasanov**
Birthday: 10.10.1955
Birthplace: Nakhchivan, Azerbaijan
Master: Mechanics, Mechanical and Mathematical Faculty, Baku State University, Baku, Azerbaijan, 1977
Doctorate: Ph.D., Mathematics, Deformed Solid Mechanics, Department of Mechanics, Institute of Mechanics and Mathematics, Baku, Azerbaijan, 1991
The Last Scientific Position: Assoc. Prof., Department of General Mathematics, Nakhichevan State University, Nakhichevan, Azerbaijan, Since 1995
Research Interests: Mechanics, Mathematics
Scientific Publications: 46 Papers, 2 Textbook, 8 Study Guide